
\input amssym.def
\input amssym

\documentstyle{amsppt}
\magnification=1200

\font\cal=cmsy10

\hoffset 0,7truecm
\voffset 1truecm
\hsize 15,2truecm
\vsize 22truecm 

\catcode`\@=11
\long\def\thanks#1\endthanks{%
  \ifx\thethanks@\empty \gdef\thethanks@{%
     \tenpoint#1} 
     \else \expandafter\gdef\expandafter\thethanks@\expandafter{%
     \thethanks@\endgraf#1}%
  \fi}
\def\keywords{\let\savedef@\keywords
  \def\keywords##1\endkeywords{\let\keywords\savedef@
  \toks@{\def\usualspace{{\it\enspace}}\tenpoint} 
  \toks@@{##1\unskip.}%
  \edef\thekeywords@{\the\toks@\frills@{%
  {\noexpand\it 
  Key words and phrases.\noexpand\enspace}}\the\toks@@}}%
  \nofrillscheck\keywords}
\def\makefootnote@#1#2{\insert\footins
 {\interlinepenalty\interfootnotelinepenalty
 \tenpoint\splittopskip\ht\strutbox\splitmaxdepth\dp\strutbox 
 \floatingpenalty\@MM\leftskip\z@\rightskip\z@
 \spaceskip\z@\xspaceskip\z@
 \leavevmode{#1}\footstrut\ignorespaces#2\unskip\lower\dp\strutbox
 \vbox to\dp\strutbox{}}}
\def\subjclass{\let\savedef@\subjclass
  \def\subjclass##1\endsubjclass{\let\subjclass\savedef@
  \toks@{\def\usualspace{{\rm\enspace}}\tenpoint}
  \toks@@{##1\unskip.}%
  \edef\thesubjclass@{\the\toks@\frills@{{%
  \noexpand\rm1991 {\noexpand\it Mathematics
  Subject Classification}.\noexpand\enspace}}%
  \the\toks@@}}%
\nofrillscheck\subjclass}
\def\abstract{\let\savedef@\abstract
 \def\abstract{\let\abstract\savedef@
  \setbox\abstractbox@\vbox\bgroup\noindent$$\vbox\bgroup
  \def\envir@end{\endabstract}\advance\hsize-2\indenti
  \def\usualspace{\enspace}\tenpoint \noindent 
  \frills@{{\smc Abstract.\enspace}}}%
 \nofrillscheck\abstract}

\outer\def\endtopmatter{\add@missing\endabstract
 \edef\next{\the\leftheadtoks}\ifx\next\empty
  \expandafter\leftheadtext\expandafter{\the\rightheadtoks}\fi
 \ifmonograph@\else
   \ifx\thesubjclass@\empty\else \makefootnote@{}%
        {\thesubjclass@}\fi
   \ifx\thekeywords@\empty\else \makefootnote@{}%
        {\thekeywords@}\fi  
   \ifx\thethanks@\empty\else \makefootnote@{}%
        {\thethanks@}\fi   
 \fi
  \pretitle
  \begingroup 
  \ifmonograph@ \topskip7pc \else \topskip4pc \fi
  \box\titlebox@
  \endgroup
  \preauthor
  \ifvoid\authorbox@\else \vskip2.5pcplus1pc\unvbox\authorbox@\fi
  \preaffil
  \ifvoid\affilbox@\else \vskip1pcplus.5pc\unvbox\affilbox@\fi
  \predate
  \ifx\thedate@\empty\else
       \vskip1pcplus.5pc\line{\hfil\thedate@\hfil}\fi
  \preabstract
  \ifvoid\abstractbox@\else
       \vskip1.5pcplus.5pc\unvbox\abstractbox@ \fi
  \ifvoid\tocbox@\else\vskip1.5pcplus.5pc\unvbox\tocbox@\fi
  \prepaper
  \vskip2pcplus1pc\relax}
\def\foliofont@{\tenrm} 
\def\headlinefont@{\tenpoint} 


\def\refstyle#1{\uppercase{%
  \if#1A\relax \def\keyformat##1{[##1]\enspace\hfil}%
  \else\if#1B\relax 
      \def\keyformat##1{\aftergroup\kern
           aftergroup-\aftergroup\refindentwd}%
      \refindentwd\parindent 
  \else\if#1C\relax 
      \def\keyformat##1{\hfil##1.\enspace}%
  \else\if#1D\relax 
      \def\keyformat##1{\hfil\llap{{##1}\enspace}} 
  \fi\fi\fi\fi}
}
\refstyle{D}   
\def\refsfont@{\tenpoint} 
\def\email{\let\savedef@\email
  \def\email##1\endemail{\let\email\savedef@
  \toks@{\def\usualspace{{\it\enspace}}\endgraf\indent\tenpoint}%
  \toks@@{\tt ##1\par}%
  \expandafter\xdef\csname email\number\addresscount@\endcsname
  {\the\toks@\frills@{{\noexpand\smc E-mail address\noexpand\/}:%
     \noexpand\enspace}\the\toks@@}}%
  \nofrillscheck\email}

\def\A{{\Bbb A}}

\def\D{\text{\rm D}}

\def\F{{\Bbb F}}

\def\K{\text{\rm K}}

\def\M{\text{\rm M}}

\def\Q{{\Bbb Q}}

\def\T{\text{\rm T}}

\def\Z{{\Bbb Z}}

\def\SL{\operatorname{SL}}
\def\GL{\operatorname{GL}}
\def\Aut{\operatorname{Aut}}
\def\End{\operatorname{End}}

\def\Char{\operatorname{char}}
\def\Ker{\operatorname{Ker}}
\def\Nor{\text{\cal N}}
\def\Nrd{\text{\rm N}}
\def\gcd{\operatorname{gcd}}
\def\proof{\noindent{\bf Proof.\ }}


\newcount\secno

\secno=1

\topmatter

\title
On the lower garland of certain\\
subgroup lattices in linear groups
\endtitle

\author
Alexandre A.$\,$Panin
\endauthor

\affil
{\it Department of Mathematics and Mechanics\\ 
St.Petersburg State University\\ 
2 Bibliotechnaya square,\\
St.Petersburg 198904, Russia}
\endaffil

\address
\endaddress
\email
alex@ap2707.spb.edu
\endemail



\thanks
\endthanks

\abstract
This paper is a revised and corrected version of [BP].
We describe here the lower garland of some lattices
of intermediate subgroups in linear groups.
The results are applied to the case of subgroup lattices
in general and special linear groups over a class of rings,
containing the group of rational points $\T$ of a maximal
non-split torus in the corresponding algebraic group.
It turns out that these garlands coincide with the
interval of the whole lattice, consisting of
subgroups between $\T$ and its normalizer.
Bibliography: 10 titles.

\endabstract

\endtopmatter

\document

\heading
Introduction
\endheading

Let $S$ be an associative ring with unit and $R$ its unitary
subring contained in the center of $S$. One can consider the
following lattice of intermediate subgroups:

$$Lat(\Aut(_SS),\Aut(_RS))=\{H:\ \Aut(_SS)\leqslant H\leqslant \Aut(_RS)\}$$
(see [BKH]). If $S$ is a free left $R$--module of rank $n$, this lattice
is isomorphic to the lattice of matrix subgroups

$$Lat(\T,\GL(n,R))=\{H:\ \T\leqslant H\leqslant \GL(n,R)\},$$
where by $\T=\T(S)$ we denote the image of $S^*$ under the
inclusion which carries each element $\alpha\in S^*$ to the matrix of
the operator of right multiplication by this element with respect to the fixed
basis of this ring extension.

\smallskip

Let $R$ be commutative and $S=R^n$. Then we can consider $R$ as a
subring of $S$ via ``diagonal'' embedding. Thus
we get the lattice of subgroups in $\GL(n,R)$,
containing the group of diagonal matrices $\D(n,R)$. The problem
of description of this lattice
has been solved for a rather wide class of rings (see the surveys [V1],
[V2] for the background information).

\smallskip

If $R$ is a field and $S$ its finite extension,
the achieved progress in the description of the lattice of intermediate
subgroups is not so notable. The cases, when $k$ is the field of
real numbers, finite or local fields are investigated (qualitative or
quantitative) completely or almost completely. The other known results
are concerned with the quadratic extensions (see also [V1], [V2]).

\smallskip

At the same time a more general problem can be stated.
Let $G'\leqslant \Aut(_RS)$. Then one can consider the
subgroup lattice

$$Lat(\Aut(_SS)\cap G',G')=
\{H':\ \Aut(_SS)\cap G'\leqslant H'\leqslant G'\}$$

Of course, one can hardly expect to obtain complete description
of this lattice for an arbitrary $G'$, but the case of finite field
extensions and classical groups $G'$ seems to be more acceptable.

\smallskip

The present paper is devoted to the investigation of certain properties
of the group $\Aut(_SS)\cap G'$. Namely, several conditions
on rings and groups which allow to calculate its normalizer in $G'$
are proposed. Under some additional restrictions the lower garland of the
lattice under consideration is described.
These results are applied to the case of general and special linear
groups.

\smallskip

Some results of this paper concerning the normalizer of 
$\Aut(_SS)\cap G'$ already appeared in [BP], but they were
extremely incomplete in the form presented there (even the
main ``conjecture'' of [BP] was trivial). Taking a chance,
the author would like to ask to be excused by a potential
reader of [BP].

\smallskip

Hereafter by a ring we mean an associative ring with unit. All subrings are
assumed to be unitary and all modules to be left ones.

\heading
\S~\the\secno. The notion of a garland
\endheading
\advance\secno by 1

Let $G$ be a group and $G_0$ its subgroup. One can be interested
in the description of the lattice of intermediate subgroups

$$Lat(G_0,G)=\{H:\ G_0\leqslant H\leqslant G\}$$
of $G$ which contain $G_0$.

Z.I.Borewicz [BKH] proposed a natural partition of this lattice
into pairwise disjoint pieces which allows to reduce the primary problem
to the problem of description of each of these pieces separately.

Namely, a structure of graph called the {\it normality graph} is introduced
on the set $Lat(G_0,G)$. We take the intermediate subgroups $H$
as the vertices of this graph and link two of them with an edge
iff one of the corresponding subgroups is normal in the other.

\smallskip

\noindent {\bf Definition.} A {\it garland} of the lattice
$Lat(G_0,G)$ is a connected component of the normality graph.

\smallskip

By this definition the lattice $Lat(G_0,G)$ is the
disjoint set-theoretical union of its garlands. Thus, to describe
the whole $Lat(G_0,G)$, it is sufficient to enumerate all the garlands and
then to investigate each garland separately.

\smallskip

Subgroups $G_0$ and $G$ are among intermediate, hence they belong to
certain garlands. The garland containing $G_0$ is called {\it lower}
and containing $G$ {\it upper} garland of the lattice $Lat(G_0,G)$.

\smallskip

One can find a number of examples of decomposition of subgroup lattices in
the union of their garlands in [AH].

\heading
\S~\the\secno. General case
\endheading
\advance\secno by 1

We shall use the following notations. Let $G$ be a group,
$G_0\leqslant G$. If $G'\leqslant G$ (an arbitrary subgroup), then
$G_0'=G_0\cap G'$.

\smallskip

Let $S$ be a ring and $R$ its subring contained in the center of $S$.
Let $M$ be an $S$--module. We denote $G=\Aut(_RM),\ G_0=\Aut(_{S}M)$.
Let $G'\leqslant G$.

\smallskip

We omit the proofs of all assertions of this section since they
can be obtained in standard way.

\proclaim
{Lemma 2.1} The following assertions are equivalent:
\smallskip\indent
$(i)$ $\Nor_{G'}G_0'\subseteq \Nor_GG_0\cap G'$;
\smallskip\indent
$(ii)$ if $G_0'\trianglelefteq H'\leqslant G'$, then
$H'\leqslant \Nor_GG_0$.
\endproclaim

\smallskip

We shall assume that to the end of this section one of these two
equivalent conditions is fulfilled.

\proclaim
{Corollary 1}
$\Nor_{G'}G_0'=\Nor_GG_0\cap G'$.
\endproclaim

\smallskip

Let's denote by $L_0$ the lattice of subgroups $Lat(G_0,\Nor_GG_0)$,
and by $L_0'$ the lattice $Lat(G_0',\Nor_{G'}G_0')$.

\proclaim 
{Corollary 2}
$L_0'=L_0\cap G'$.
\endproclaim

\smallskip

Let now $S$ be a commutative ring, $M=S$. Then the group
$G_0\thickapprox S^{*}$ is Abelian. 

\proclaim 
{Corollary 3}
If $G_0$ is a maximal Abelian subgroup in $G$
(e.g. if $S$ is additively generated by its invertible elements),
then $G_0'$ is also a maximal Abelian subgroup in $G'$.
\endproclaim

\vfill\eject

\heading
\S~\the\secno. Calculation of the normalizer
\endheading
\advance\secno by 1

Note that $\End(_SM)$ is an $R$--module since $R$ is contained in the
center of $S$.
We assume that the following property holds true:

\bigbreak

$\quad\quad\quad G_0' {\hbox {\it \ additively generates (over }} R)\
{\hbox {\it the group\ }} G_0$ 
\hfill $(+)$

\bigbreak

\proclaim
{Lemma 3.1}
$\Nor_{G'}G_0'=\Nor_GG_0\cap G'$.
\endproclaim
\proof Let $h\in \Nor_{G'}G_0',
\ t\in G_0$. Then $t=\sum\limits_{i=1}^s\alpha_it_i$, where
$\alpha_i\in R,\ t_i\in G_0'$. Thus
$h^{-1}th=\sum\limits_{i=1}^s\alpha_ih^{-1}t_ih=
\sum\limits_{i=1}^s\alpha_i^{}t_i'$ for some $t_i'\in G_0'$.
Hence $h^{-1}th\in G_0$. The opposite inclusion is trivial.

\proclaim
{Theorem 3.2}
Let $S$ be a ring and $R$ its subring contained
in the center of $S$. Assume that $S$ is additively generated
by its invertible elements. Let $M$ be a free $S$--module
of finite rank, $G'$ be a subgroup of $G$ such that the condition
$(+)$ is fulfilled. Then the normalizer of $G_0'$ in $G'$ is
equal to the intersection of the semidirect product of the
normal subgroup $G_0$ and the subgroup isomorphic to the group $\Aut(S/R)$ of
all ring automorphisms of $S$, identical on $R$, with
the group $G'$. 
\endproclaim
\proof This theorem was proved in the paper of V.A.Koibaev [K] for $G'=G$.
Therefore it is sufficient to apply Lemma 3.1.

\smallskip

We are primarily interested in the case $M=S$ (see Introduction).
But the examination of the condition $(+)$ even in this case is not possible
for a sufficiently large class of rings. Moreover, this condition
does not allow to obtain further information on the structure
of the lower garland of the lattice $Lat(G_0',G')$. So we shall try
to impose another restriction on the class of rings under consideration.

\smallskip

Let $S$ be a ring and $R$ its subring, which is an integral domain
contained in the center of $S$. Let also $S$ be a free
$R$--module of finite rank with a basis $\omega_1,\ldots,\omega_n$.

We consider an embedding of $S$ into the ring of matrices $\M(n,R)$:
each element $\alpha=\alpha_1\omega_1+\ldots +\alpha_n\omega_n\in S$
(where $\alpha_i\in R$) maps to a matrix $t(\alpha)=(t_{ij}(\alpha))$,
where $\omega_i\alpha=\sum\limits_{j=1}^nt_{ji}(\alpha)\omega_j$.

Let's denote $\T=t(S^*)$. Let $k$ be
the quotient field of $R$ and $\overline{k}$ its algebraic closure.
We spread the mapping $t$ to an affine space $\A^n$, changing
$R$ to $\overline{k}$ and $S$ to $\overline{k}\otimes_RS$. 
It follows from the definition of the mapping $t$ that
$t(\alpha\beta)=t(\beta)t(\alpha)$ for every
$\alpha,\beta\in\overline{k}\otimes_RS$.
We put $\widetilde \T=t((\overline{k}\otimes_RS)^*)$. It's clear that
$\widetilde \T$ is a group.

\proclaim
{Lemma 3.3}
$\widetilde \T$ is a connected $k$--defined algebraic group.
\endproclaim

\proof It's clear that $\widetilde \T=t(U)$,
where $U$ is a Zariski open subset of $\A^n$.
Therefore $\widetilde T$ is a connected algebraic group.

The set $U$ is not changed after passing to $S^{op}$, and we get
a morphism $t\vert _U:U\to \GL(n,\overline{k})$ of $k$--groups, defined
over $k$. Then its image is a $k$--subgroup of $\GL(n,\overline{k})$.

\smallskip

Let $G'$ be an arbitrary subgroup of $\GL(n,R)$.
We put $\T'=\T\cap G'$. Assuming that $\overline{k}$ is naturally
embedded in $\overline{k}\otimes_RS$, we impose the following
restriction on $G'$: 

\bigbreak

$\overline{\T}$ {\it is additively generated over} $\overline{k}$
{\it by the connected component of} $\overline{\T'}$ \hfill $(*)$

\proclaim
{Lemma 3.4}
$\Nor_{G'}\T'=\Nor_{\GL(n,R)}\T\cap G'$.
\endproclaim

\proof Let $g\in \Nor_{G'}\T'$, i.e.
$g\T'g^{-1}=\T'$. Since the conjugation is an isomorphism of algebraic
groups, we see that $g\overline{\T'}g^{-1}=\overline{\T'}$. By the
condition $(*)$
$g\T g^{-1}\subseteq g\overline{\T}g^{-1}\subseteq\widetilde \T$.
Passing to ``R--points'', we get $g\in \Nor_{\GL(n,R)}\T$.

\smallskip

\noindent {\bf Remark.} In the course of the proof of Lemma 3.4 we used only
the fact that $\overline{\T}$ is additively generated over $\overline{k}$
by the whole group $\overline{\T'}$.

\smallskip

Now we can formulate the analogue of Theorem 3.2.

\proclaim
{Theorem 3.5}
Let $S$ be a ring and $R$ its subring, which is an integral domain
contained in the center of $S$. Let $S$ be a free 
$R$--module of finite rank which is additively generated by its
invertible elements. Let $G'$ be a subgroup of $\GL(n,R)$ such
that the condition $(*)$ is fulfilled. Then the normalizer of $\T'$
in $G'$ is equal to the intersection of the semidirect product of the
normal subgroup $\T$ and the group $\Aut(S/R)$ of all ring automorphisms
of $S$, identical on $R$, with the group $G'$.
\endproclaim

\heading
\S~\the\secno. Calculation of the lower garland
\endheading
\advance\secno by 1

Let, as above, $S$ be a ring and $R$ its subring, which is an integral domain
contained in the center of $S$. We assume that $S$ is a free
$R$--module of rank $n$.

Let $G'=\widetilde G(R)$, where $\widetilde G$ is a closed subgroup
of $\GL(n,\overline{k})$.

\proclaim
{Theorem 4.1}
Suppose that
\smallskip\indent
$(i)$ $\widetilde G$ is a connected reductive group;
\smallskip\indent
$(ii)$ $\widetilde \T\cap\widetilde G$ is a maximal torus in $\widetilde G$;
\smallskip\indent
$(iii)$ the group $\T'$ is Zariski dense in $\widetilde \T\cap\widetilde G$.
\smallskip\noindent
Then the lower garland of the lattice $Lat(\T',G')$ coincides with
the interval $Lat(\T',\break\Nor_{G'}\T')$.
\endproclaim
\proof One has to check that
$\T'\leqslant H\leqslant \Nor_{G'}\T'$ implies
$\Nor_{G'}H\leqslant \Nor_{G'}\T'$.

Let $g\in \Nor_{G'}H$. Then
$g\T'g^{-1}\subseteq H\subseteq \Nor_{G'}\T'$.
Since the conjugation is an isomorphism of algebraic groups, we obtain
$\Nor_{G'}\T'\subseteq \Nor_{\widetilde G}\overline{\T'}$.
The latter group is closed, hence
$g\overline{\T'}g^{-1}\subseteq \Nor_{\widetilde G}\overline{\T'}$.
Since $\T'$ is dense in $\widetilde \T\cap\widetilde G$, we get
$g(\widetilde \T\cap\widetilde G)g^{-1}\subseteq 
\Nor_{\widetilde G}(\widetilde \T\cap\widetilde G)$.
It follows from $(ii)$ that the (regular) torus $\widetilde \T\cap\widetilde G$
coincides with its centralizer in $\widetilde G$. But then
the connected component
$(\Nor_{\widetilde G}(\widetilde \T\cap\widetilde G))^o$ is
equal to $\widetilde \T\cap\widetilde G$. Since any morphism of algebraic
groups maps the connected component of the first group onto the connected
component of the second, we see that
$g\in \Nor_{\widetilde G}(\widetilde \T\cap\widetilde G)$.
Passing to the ``$R$--points'', we get $g\in \Nor_{G'}\T'$.

\smallskip

Note that in the course of the proof of Theorem 4.1 we did not use the explicit
form of the normalizer of $\T'$. Instead of this we imposed
several restrictions on the groups, particularly, we assumed
$\widetilde \T\cap\widetilde G$ to be a maximal torus in $\widetilde G$.
Now we prove an analogue of Theorem 4.1 under some other conditions,
and this allows us to get some corollaries which can not be
deduced from Theorem 4.1.

\proclaim
{Theorem 4.2}
Assume that $S$ is additively generated by its
invertible elements and the group $\Aut(S/R)$ consisting of the ring
automorphisms of $S$, which are identical on $R$, is finite. If the
condition $(*)$ for $G'$ is fulfilled, then 
the lower garland of the lattice $Lat(\T',G')$ coincides with
the interval $Lat(\T',\Nor_{G'}\T')$.
\endproclaim
\proof The crucial idea in the proof of Theorem 4.1 lies in the fact that
the connected component of the normalizer of $\overline{\T'}$ coincides
with $\overline{\T'}$. Now the group
$\overline{\T'}$ may be not connected, but the scheme of the
proof is the same.

Let $g\in \Nor_{G'}H$. It follows from Theorem 3.5 that
$\Nor_{G'}\T'={\bigcup\limits_{i=1}^s}^{\varnothing}\T'A_i$
for some $A_i\in \GL(n,R)$ and $A_1=1$. Thus
$g\overline{\T'}g^{-1}\subseteq \bigcup\limits_{i=1}^s\overline{\T'}A_i$.
The latter set is a group, moreover it's a (closed) subgroup of
$\GL(n,\overline{k})$ and the index of its subgroup $(\overline{\T'})^o$ is
finite, hence this subgroup coincides with its connected component. Thus,
as it was mentioned above,
$g(\overline{\T'})^og^{-1}\subseteq (\overline{\T'})^o$.
Using $(*)$, we get $g\T'g^{-1}\subseteq \widetilde \T\cap
\widetilde G$. Passing to the ``$R$--points'', we have
$g\in \Nor_{G'}\T'$.

\proclaim
{Corollary}
In the settings of Theorem 4.1 or Theorem 4.2
$\ \Nor_{G'}\Nor_{G'}\T'=\Nor_{G'}\T'$.
\endproclaim

\noindent{\bf Remark.} Usually these two theorems
do not give any information on the structure of the lower garland
in the case of noncommutative ring $S$. Indeed, then
$\widetilde \T\cap \widetilde G$ is not a torus, and
the group $\Aut(S/R)$ contains a lot of inner automorphisms,
hence is not finite.
 
\heading
\S~\the\secno. Separable algebras
\endheading
\advance\secno by 1

Here we state without proofs some useful results on separable
algebras. See [CR] for a more detailed exposition.

\smallskip

Let $k$ be a field. All $k$--algebras assumed to be finite-dimensional.

\smallskip

\noindent{\bf Definition}. A $k$--algebra $A$ is called separable
(over $k$), if the $K$--algebra $K\otimes_kA$ is semisimple for any extension
$K$ of the field $k$. In particular, any separable algebra is automatically
semisimple.

\proclaim
{Lemma 5.1}
A $k$--algebra $A$ is separable iff 
the $\overline{k}$--algebra $\overline{k}\otimes_kA$ is semisimple.
\endproclaim

\proclaim
{Theorem 5.2}
Let $A=A_1\oplus\ldots\oplus A_s$ be a decomposition
of a semisimple $k$--algebra $A$ into the direct sum of simple algebras
$\{A_i\}$ with the centers $\{C_i\}$. Then $A$ is separable iff
all $C_i$ are separable extensions of $k$.
\endproclaim

\heading
\S~\the\secno. The case of the general linear group
\endheading
\advance\secno by 1

\proclaim
{Lemma 6.1}
Let $S$ be a commutative ring and $R$ its subring,
which is an integral domain. Assume that $S$ is a free $R$--module
of rank $n$. If the $\overline{k}$--algebra $\overline{k}\otimes_RS$
is semisimple, then $\widetilde \T$ is a maximal torus in $\GL(n,\overline{k})$
defined over $k$.
\endproclaim
\proof In view of Lemma 3.3 and dimension arguments
it remains to show that the group $\widetilde \T$ is 
diagonalizable.

It follows from the settings that $\overline{k}\otimes_RS$ is isomorphic
to the direct sum of $n$ copies of the field $\overline{k}$. Thus
$\widetilde \T$ is conjugated to the group of diagonal matrices
$\D(n,\overline{k})$.

\smallskip

\noindent{\bf Remark.} If the $\overline{k}$--algebra $\overline{k}\otimes_RS$
contains nilpotents, then $\widetilde \T$ has unipotent elements,
hence it is not a torus.

\proclaim
{Theorem 6.2}
Let $S$ be a commutative ring and $R$ its subring,
which is an integral domain. Assume that $S$ is a free
$R$--module of rank $n$, the $\overline{k}$--algebra
$\overline{k}\otimes_RS$
is semisimple (where $k$ is the quotient field of $R$), and
$\T$ is Zariski dense in $\widetilde \T$. Then the lower
garland of the lattice $Lat(T,\GL(n,R))$ coincides with the interval
$Lat(\T,\Nor_{\GL(n,R)}\T)$.
\endproclaim
\proof Follows from Theorem 4.1.

\smallskip

\noindent {\bf Remark.} One can use also Theorem 4.2: it is sufficient
to require that $S$ is additively generated by its invertible elements,
and the $k$--algebra $k\otimes_RS$ is semisimple, but then
$\widetilde \T$ may be not a torus.

\proclaim
{Corollary}
Let $R=k$ be a field, $S=K_1\oplus\ldots\oplus K_t$,
where $K_i/k$ are finite extensions of $k$.
Then:
\smallskip\indent
$(i)$ if all $K_i/k$ are separable, then
$\widetilde \T$ is a maximal torus in $\GL(n,\overline{k})$
defined over $k$;
\smallskip\indent
$(ii)$ if $k$ is infinite, then the lower garland of the lattice $Lat(\T,
\GL(n,k))$ coincides with the interval $Lat(\T,\Nor_{\GL(n,k)}\T)$.
\endproclaim
\proof Follows from Lemma 6.1 and results of \S$\,$5.

\smallskip

The case of a finite field will be treated in \S$\,$8.

\smallskip

Note that this Corollary generalizes the results of Al Hamad [AH].

\smallskip

\noindent {\bf Example}. Let $R=\Z$ be the ring of integer numbers.
Consider the ring $S=\Z\oplus \Z w$, where $w^2=0$. It's clear
that there are only two automorphisms of $S/R$, namely, they map
$w$ to $\pm w$. Furthermore, the invertible elements of $S$ are of the
form $\pm 1+bw$, where $b\in R$.
Thus the first two conditions of Theorem 4.2 are fulfilled.

The group $\T$ is equal (if we choose elements $1,w$ as a basis of $S/R$) to
the set of matrices

$$\left\{ {\pmatrix a&0\\ b&a\\ \endpmatrix : a=\pm1, b\in R}\right\}$$

To get $\overline{\T}$, one has to take an arbitrary element
of $\overline{\Q}$ in place of $b$ in the latter expression. We note that
this group is not connected, but the condition $(*)$ is nevertheless
fulfilled. Thus we can apply Theorem 4.2.

Note that it is impossible to apply Theorem 4.1 in this situation
because of the existence of the unipotent elements in $\T$.

\heading
\S~\the\secno. The case of the special linear group
\endheading
\advance\secno by 1

Let $\widetilde G=\SL(n,\overline{k})$.

\proclaim
{Lemma 7.1}
Let $S$ be a commutative ring and $R$ its subring
which is an integral domain. Let $S$ be a free $R$--module
of rank $n$. If the $\overline{k}$--algebra $\overline{k}\otimes_RS$ is
semisimple, then $\widetilde \T\cap\widetilde G$ is a
maximal $k$--defined torus in $\SL(n,\overline{k})$.
\endproclaim
\proof We show that the group $\widetilde \T\cap\widetilde G$ is connected.

Since $\overline{k}\otimes_RS$ is semisimple, we see that
$\widetilde \T\cap\widetilde G$ is conjugated with the subgroup of
$\D(n,\overline{k})\cap\widetilde G$ consisting of diagonal matrices
with the determinant $1$. The latter group is connected since
a polynomial $x_1x_2\ldots x_n-1$ is irreducible over $\overline{k}$.

To complete the proof it is sufficient to turn to the dimension
arguments.

\proclaim
{Theorem 7.2}
Let $S$ be a ring and $R$ its subring which
is an integral domain containing in the center of $S$. Assume that $S$
is a free $R$--module of rank $n$ which is additively generated
by its invertible elements, and $\T'$ is Zariski dense in
$\widetilde \T\cap\widetilde G$. Then the normalizer of the subgroup $\T'$
in $G'$ is equal to the intersection of the semidirect product
of the normal subgroup $\T$ and the group $\Aut(S/R)$ of the ring automorphisms
of $S$, identical on $R$, with $G'$.
\endproclaim
\proof It is clear that each matrix of
$\widetilde \T$ can be transformed to an element of
$\widetilde \T\cap\widetilde G$ via the multiplication by a scalar matrix,
therefore the condition $(*)$ from \S$\,$3 is verified and we can apply
Theorem 3.5 (see the remark after Lemma 3.4).

\proclaim
{Theorem 7.3}
Let $S$ be a commutative ring and $R$ its subring which is an
integral domain. Assume that $S$ is a free $R$--module of rank $n$,
the $\overline{k}$--algebra $\overline{k}\otimes_RS$ (where
$k$ is the quotient field of $R$) is semisimple, and
$\T'$ is Zariski dense in $\widetilde \T\cap\widetilde G$.
Then the lower garland of the lattice $Lat(\T',\SL(n,R))$ is equal to
the interval $Lat(\T',\Nor_{\SL(n,R)}\T')$.
\endproclaim

\proclaim
{Corollary}
Let $R=k$ be a field, $S=K_1\oplus\ldots\oplus K_t$,
where $K_i/k$ are finite separable extensions of $k$.
Then:
\smallskip\indent
$(i)$ $\widetilde \T\cap\widetilde G$ is a maximal $k$--defined torus
in $\SL(n,\overline{k})$;
\smallskip\indent
$(ii)$ if $k$ is infinite, then the lower garland of the lattice $Lat(\T',
\SL(n,k))$ is equal to the interval $Lat(\T,\Nor_{\SL(n,k)}\T')$,
which also equals the intersection of the lower garland of the lattice
$Lat(\T,\GL(n,k))$ with $\SL(n,k)$.
\endproclaim

\smallskip
The case of a finite field will be treated in \S$\,$8.

\smallskip

Now we analyse several examples.

\smallskip

$1^0$. If $K/k$ is a finite field extension, then $\T'=\Ker\Nrd_{K/k}$.

\smallskip

$2^0$. Let $k=\Q\/, K=\Q(\sqrt d)$ be its quadratic extension,
where $d$ is a squarefree integer rational number. If we take
$1,\sqrt d$ as a basis of $K/k$, then

$$\T'=\left\{ {\pmatrix x&yd\\ y&x\\ \endpmatrix :\ x,y\in\Q,
\ x^2-dy^2=1} \right\}$$

It is clear that there are only two automorpisms of $S/R$, namely,
identical one and the linear mapping which carries $\sqrt d$ to $-\sqrt d$.

Now we can explicitly calculate the subgroup
$\Nor_{\SL(2,\Q)}\T'$. Indeed, it follows
from the Corollary of Theorem 7.3 that the result depends on
the arithmetical properties of $d$ (namely, on the solvability
in rational numbers of the equation $x^2-dy^2=-1$):
if $d$ is positive and do not have any prime divisors of the form
$4m+3$, then

$$\Nor_{\SL(2,\Q)}\T'=
\T'\cup \left\{\pmatrix x&-yd\\ y&-x\\ \endpmatrix:\ x,y\in\Q,
\ x^2-dy^2=-1 \right\}=$$

$$=\T'\cup \T'\pmatrix x_0&-y_0d\\ y_0&-x_0\\ \endpmatrix,$$
where $(x_0,y_0)$ is a fixed solution of the equation $x^2-dy^2=-1$;
otherwise $\Nor_{\SL(2,\Q)}\T'\break=\T'$.
Thus, the lower garland of the lattice $Lat(\T',\SL(2,\Q))$
consists either of one or of two subgroups.

\smallskip

$3^0$. Let $K/k$ be a pure inseparable field extension of degree
$q=p^n$, where $\Char k=p$. Since $\Nrd_{K/k}(\alpha)=\alpha^q$, we see that
$\Ker\Nrd_{K/k}=1$. Indeed, $\alpha^q=1$ iff $(\alpha-1)^q=0$.

\heading
\S~\the\secno. An elementary approach
\endheading
\advance\secno by 1

Let $K$ be a finite separable extension of an infinite field $k$,
$N$ a fixed natural number.

\proclaim
{Lemma 8.1}
Let $x\in K\backslash k$. Then the set of 
$\alpha\in k$ such that $(x+\alpha)^N\in k$ is finite, moreover,
it can contain not more than $N$ elements.
\endproclaim

\proof We denote by $p$ the characteristic exponent of $k$.
Then $N=p^rN'$, where $\gcd(N',p)=1$. In view of the
separability of $K/k$ we can assume that $\gcd(N,p)=1$.

Let's suppose that the assertion of Lemma is false and choose
different $\alpha_0,\ldots,\break\alpha_N\in k$ such that
$(x+\alpha_i)^N\in k$ for every $i=0,1,\ldots,N$. Using the
standard binomial formula, we get a system of equations
of the form $CX=Y$, where $C\in \GL(N+1,k),\ Y\in k^{N+1}$,
whence $X\in k^{N+1}$, and, in particular, $Nx\in k$, therefore $x\in k$.
A contradiction.

\proclaim
{Corollary}
There exist infinitely many primitive elements
of the extension $K/k$ with the norm $1$.
\endproclaim

\proof Let's take a primitive element $x$ of
the finite separable extension $K/k$. For every $\alpha\in k$
one can consider the element
$x_\alpha={(x+\alpha)^n\over \Nrd_{K/k}(x+\alpha)}$, where $n=(K:k)$.
All such elements have the norm $1$. Since there are only finitely
many intermediate subfields
$F,\ k\subseteq F\subseteq K$, it follows from Lemma 8.1
that there are infinitely many different primitive among the
elements $x_\alpha$.

\smallskip

If $k$ is a finite field, then, as it can be easily verified,
there exists a primitive element of the extension $K/k$ with the norm $1$.

\smallskip

Let $R=k$ be a finite field, $S=K_1\oplus\ldots\oplus K_t$, where $K_i/k$
are finite field extensions of degree $n_i$,
$\widetilde G=\SL(n,\overline{k})$.

\proclaim
{Lemma 8.2}
If $S\neq \F_3\oplus{\F}_3$, then
$\Nor_{G'}\T'=\Nor_{\GL(n,k)}\T\cap G'$.
\endproclaim

\proof We show that the set $S'$ which consists of
$k$--linear combinations of the elements $s\in S$ such that $t(s)\in G'$
contains $S^*$. Then the condition $(*)$ is fulfilled, and we can use
Lemma 3.4.

A. Let not all of $K_i/k$ be $1$--dimensional. For every index $j$ such
that $n_j\geqslant 2$ we choose a primitive element $\alpha_j$ of
the extension $K_j/k$ with the norm $1$. Then
$(0,\ldots,\alpha_j-1,\ldots,0)\in S'$, and we get $K_j\subseteq S'$.

If some of $K_i$ are $1$--dimensional, then there is nothing to prove for
$k=\F_2$. If $k\neq \F_2$, then for every $j$ such that
$n_j=1$, one can choose an element $c_j\in K_j=k$ not equal to $0$ and $1$.
There is no loss in assuming $n_1\geqslant 2$. In view of the
surjectivity of the norm homomorphism one can find an element
$\alpha\in K_1$ with the norm $c^{-1}$.
Then it is clear that $(0,\ldots,c_j-1,\ldots,0)\in S'$, therefore
$K_j\subseteq S'$.

B. Let all $K_i=k$ and $\vert k\vert\geqslant 4$. Let's choose $c\in k^*$
such that $c^2\neq 1$. Then $(c^2-1,c-1,\ldots,c-1,0)\in S'$,
whence $(c+1,1,\ldots,1,0)\in S'$. If we act in the same manner
with $c+1\neq 0$, we get $(1,0,\ldots,0)\in S'$.

C. It remains to examine the case $K_i=\F_3,\ t\geqslant 3$.
If $t$ is odd, then $(1,-1,\ldots,\break -1)+(1,\ldots,1)=(-1,0,\ldots,0)\in
S'$.
If $t$ is even, then by the same reasons $(-1,-1,0,\ldots,0)\in S'$, and also
$(-1,1,0,\ldots,0)=(1,-1,1,-1,\ldots,-1)+(1,-1,-1,1,\ldots,1)\in S'$.
Lemma is proved.

\smallskip

\noindent{\bf Remark.} It's easy to verify that for
$S=\F_3\oplus \F_3$ the assertion of Lemma 8.2 doesn't hold true.

\proclaim
{Corollary}
Let $R=k$ be a field, $S=K_1\oplus\ldots\oplus K_t$,
where $K_i/k$ are finite separable field extensions, and:
\smallskip\indent
$(i)$ not more than two of them are equal to $\F_2$;
\smallskip\indent
$(ii)$ $S\neq \F_3\oplus \F_3$.
\smallskip\indent
Then the normalizer of\/ $\T'$ in $G'$ is equal to
the intersection of the semidirect product of $\T$ and $\Aut(S/R)$
with the group $G'$. 
\endproclaim
\smallskip

G.Seitz [S] obtained an almost exhaustive description of
the intermediate subgroups in the finite groups of Lie type
which contain a group of rational points of a maximal
torus in the corresponding algebraic group. For simplicity we state his
results only for Chevalley groups.

Let $\widetilde G=G(\Phi,\overline{k})$ be a Chevalley group of type $\Phi$
over the algebraic closure $\overline{k}$ of a finite field
$k=\F_q$, where $\Char k\neq 2,3$ and $q\geqslant 13$, $G=G(\Phi,k)$.
It was proved in [S] that the group of $k$--rational points $\T$
of an arbitrary maximal torus $\widetilde \T$ in $\widetilde G$ is
a paranormal subgroup of $G$ (see [BB], [V2]), which means, in particular, that
the lower garland of the lattice of intermediate subgroups of $G$
containing $\T$ coincides with the interval $(\T,\Nor_G\T)$.

Therefore one can elaborate the statements of the corresponding
corollaries from Theorems 6.2 and 7.3.

\smallskip

For those fields $k$ which do not satisfy the conditions under which the
results of Seitz hold true, the lower garland of the lattice
$Lat(\T,G)$ actually may not coincide with the interval $Lat(\T,\Nor_G\T)$.
The examples of such kind can be easily constructed
for $k=\F_3$, $K=(\F_3)^n$, $n\geqslant 2$.

\heading
Acknowledgements
\endheading

This work has been carried out in the framework
of the {\it ``R\'eseau de formation--recherche'' in $\K$--theory
and linear algebraic groups} during the author's visit
at the Universit\'e de Franche-Comt\'e, Besan\c con, France.
The author is personally grateful to E.Bayer--Fluckiger for the
hospitality.

Also the author would
like to express his thanks to B.Kunyavskii and V.Chernou\-sov
for numerous discussions, N.Vavilov for the thorough reading
of the manuscript of this paper, and ``The Rolling Stones''
for the wonderful music which was the main source of
inspiration during the preparation of this paper.

\enddocument